\documentclass{ip-journal}
\usepackage{amsmath}
\usepackage{amssymb}
\usepackage{amsthm}
\usepackage{mathrsfs}
\usepackage{pdfpages}
\usepackage{dsfont}
\usepackage{lineno}
\usepackage{enumerate}
\usepackage{graphicx}
\usepackage{epstopdf}
\usepackage{color}
\usepackage{xspace}
\usepackage{subfig}
\usepackage{pdfsync}
\usepackage{dsfont}

\DeclareGraphicsExtensions{.pdf}
\graphicspath{{./Fig/}}

\newcommand{\bq}{\begin{equation}}
\newcommand{\eq}{\end{equation}}
\newcommand{\R}{\mathbb{R}}

\newcommand{\abs}[1]{\left\vert#1\right\vert}

\newcommand{\MA}{Monge-Amp\`ere\xspace}

\newcommand*\Laplace{\mathop{}\!\mathbin\bigtriangleup}

\newtheorem{theorem}{Theorem}
\theoremstyle{lemma}

\newtheorem{definition}{Definition}

\theoremstyle{remark}

\makeindex

\begin{document}
\title{Seismic Inversion and the Data Normalization for Optimal Transport}

\author{bj\"orn engquist}
\address{Department of Mathematics and ICES, The University of Texas at Austin, 1 University Station C1200, Austin, TX 78712 USA}
\email{engquist@math.utexas.edu}

\author{yunan yang}

\thanks{
It is our honor to dedicate this paper to Professor Roland Glowinski on the occasion of his eighties birthday.
We thank Junzhe Sun and Lingyun Qiu for constructive discussions and thank the sponsors of the Texas Consortium for Computational Seismology (TCCS) for financial support. The first author was partially supported by NSF DMS-1620396.
} 
\address{Department of Mathematics, The University of Texas at Austin, 1 University Station C1200, Austin, TX 78712 USA}
\email{yunanyang@math.utexas.edu}

\begin{abstract}
Full waveform inversion (FWI) has recently become a favorite technique for the inverse problem of finding properties in the earth from measurements of vibrations of seismic waves on the surface. Mathematically, FWI is PDE constrained optimization where model parameters in a wave equation are adjusted such that the misfit between the computed and the measured dataset is minimized. In a sequence of papers, we have shown that the quadratic Wasserstein distance from optimal transport is to prefer as misfit functional over the standard $L^2$ norm. Datasets need however first to be normalized since seismic signals do not satisfy the requirements of optimal transport. There has been a puzzling contradiction in the results. Normalization methods that satisfy theorems pointing to ideal properties for FWI have not performed well in practical computations, and other scaling methods that do not satisfy these theorems have performed much better in practice. In this paper, we will shed light on this issue and resolve this contradiction.
\end{abstract}

\date{January 6, 2018}
\maketitle
%\tableofcontents

\section{Introduction}
%1. Introduction (Seismic inversion and briefly describe results)

There are two major processes in exploration seismology. One is migration or reverse time migration (RTM), which determines details of the reflecting surfaces assuming an approximate model of wave velocity~\cite{Baysal1983}. Seismic inversion or full waveform inversion (FWI) is a process of recovering the quantitative features of the geophysical structure. The focus is currently on the nonlinear inverse problem of building an accurate model of the wave velocity in the earth. This is done in an iterative process where a forward seismic simulation based on the unknown velocity is matched to the actual recordings~\cite{Mora1989}. There are many related techniques in seismic exploration. Wave equation travel time tomography~\cite{luo1991wave} and the ray-based tomography are phase-like inversion methods~\cite{Schuster2011}. Least-squares inversion is known as linearized waveform inversion~\cite{lailly1984migration,tarantola1984linearized} and the least-square reverse time migration (LSRTM)~\cite{Dai2012} based on the Born approximation~\cite{hudson1981use, van1954correlations} is one example, where the background model is not updated after each iteration.

FWI is a high-resolution seismic imaging technique, which recently has been getting great attention from both academia and industry~\cite{Virieux2017}. The goal of FWI is to find both the small-scale and the large-scale components, which describe the geophysical properties using the entire content of seismic traces. A trace is the time history of seismic vibrations measured at a receiver. In this paper, we will consider the inverse problem of finding the wave velocity of an acoustic wave equation in the interior of a domain from knowing the Cauchy boundary data together with natural boundary conditions~\cite{Clayton1977}. This is implemented by minimizing the difference or mismatch between computed and measured data on the boundary. It is thus a partial differential equation (PDE) constrained optimization.

FWI is increasing in popularity even if it is still facing major computational challenges. Depending on the parameterization of the velocity model this inverse PDE-constrained optimization problem is often highly non-unique and non-convex in nature. The least-squares norm ($L^2$), which is classically used in FWI to measure the data mismatch, suffers from local minima trapping, the so-called cycle skipping issues, and sensitivity to noise~\cite{Seismology2011}. We will see that optimal transport based quadratic Wasserstein metric ($W_2$) is capable of dealing with some of these limitations by including both amplitudes mismatches and travel time differences.

%
%There are various kinds of numerical techniques that are used in seismic inversion, but FWI is increasing in popularity even if it is still facing major computational challenges. Also as PDE-constrained optimization, the problem is highly non-unique and non-convex in nature. The least-squares norm ($L^2$), classically used in FWI, suffers from local minima trapping, the so-called cycle skipping issues, and sensitivity to noise~\cite{Seismology2011}. We will see that optimal transport based quadratic Wasserstein metric ($W_2$) is capable of dealing with the last two limitations by including both amplitudes mismatches and travel time differences.

The idea of using Wasserstein metric for seismic inversion was first proposed in \cite{EFWass}. This metric is based on optimal transport~\cite{Villani}. We first transform our datasets of seismic signals into density functions of two probability distributions. Next, we find the optimal map between these two datasets and compute the corresponding transport cost as the misfit function in FWI, either by solving a \MA equation for the entire dataset or by using the explicit 1D formula~\cite{Villani} measuring the misfit trace by trace~\cite{yang2017application}. Following the idea that changes in velocity cause a shift or ``transport'' in the arrival time of a seismic signal, we demonstrated in~\cite{engquist2016optimal} the advantageous mathematical properties of the quadratic Wasserstein metric ($W_2$) and provided rigorous proofs that laid a solid theoretical foundation for this new misfit function. 

%
%
%The idea of using optimal transport for seismic inversion was first proposed in~\cite{EFWass}. The Wasserstein metric is a concept based on optimal transportation~\cite{Villani}. Here, we transform our datasets of seismic signals into density functions of two probability distributions. Next, we find the optimal map between these two datasets and compute the corresponding transport cost as the misfit function in FWI by either solving the \MA equation or using the explicit 1D formula~\cite{yang2017application}.  

%Following the idea that changes in velocity cause a shift or ``transport''  in the arrival time, \cite{engquist2016optimal} demonstrated the advantageous mathematical properties of the quadratic Wasserstein metric ($W_2$) and provided rigorous proofs that laid a solid theoretical foundation for this new misfit function. We can apply $W_2$ as misfit function in two different ways: trace-by-trace comparison which is related to 1D optimal transport, and the entire dataset comparison in multiple dimensions. We will see that solving the \MA equation in each iteration of FWI is a useful technique~\cite{yang2017application}. An analysis of the 1D optimal transport approach and the conventional misfit functions such as $L^2$ norm and integral $L^2$ norm illustrated the intrinsic advantages of this transport idea~\cite{yangletter}. 

There are two main requirements for signals $f$ and $g$ in optimal transport theory:
\bq
f(t)\geq 0,\ g(t)\geq 0,\ <f> = \int f(t) dt = \int g(t) dt = <g>.
\eq
Since these constraints are not expected for seismic signals, some data pre-processing is needed before we can implement the Wasserstein-based FWI. In~\cite{yangletter,yang2017application} we normalized the signals by adding a constant,
\bq\label{eq:linear}
\tilde{f}(t) = \frac{f(t) + c}{<f+c>},\ \tilde{g}(t) = \frac{g(t) + c}{<g+c>},\ c = \min_t(f(t),g(t)).
\eq
This worked remarkably well in realistic large scale examples~\cite{yangletter,yang2017application} together with the adjoint-state method for optimization in either the 1D or the \MA based techniques. This linear normalization does, however, not give a convex misfit functional with respect to simple shifts. Other normalizations that generate convex misfits were also tried as, for example, only using the positive part of the signals, squaring or taking the envelope or the absolute values~\cite{EFWass,engquist2016optimal}. It was puzzling that these misfit functionals performed poorly with realistic datasets.

FWI will be introduced in section two, and we will present relevant parts of optimal transport theory as background in section three. The new material is in section four where data normalizations are discussed. We will see that it is desirable to require the scaling function to be differentiable so that it is easy to apply chain rule when calculating the  Fr\'{e}chet derivative for FWI backpropagation and also better suited for the \MA solver. Other aspects of normalization are also discussed that explain the contradictions mentioned above and finally ending up with a new normalization that satisfies most of the essential properties:
\begin{equation}\label{eq:mixed}
     \tilde{f}(t) = \left\{
     \begin{array}{rl}
     &  (f(t) + \frac{1}{c}) / b,\ f(t) \geq 0,\ c>0   \\
     & \frac{1}{c} \exp(cf(t)) / b,\ f(t)<0  \end{array} \right.
\end{equation}
where $b = <(f + \frac{1}{c})\mathds{1}_{f\geq 0} + \frac{1}{c} \exp(cf) \mathds{1}_{f< 0}  >$.

%Since these are not expected for seismic signals, some data pre-processing is needed before we can implement Wasserstein-based FWI. A natural way is to separate the data into positive and negative parts. This approach preserves the desirable theoretical properties of convexity to shifts and noise insensitivity, but it is not easy to combine with the adjoint-state method~\cite{Plessix} and more realistic examples. We require the scaling function to be differentiable so that it is easy to apply chain rule when calculating the Fr\'{e}chet derivative for FWI backpropagation and also better suited for the \MA solver. There are some differentiable scaling functions to rescale the datasets so that they become positive which also preserve the convexity concerning simple shifts, but they create more non-uniqueness and more local minima are observed in inversion than normal. The convexity of $W_2$ highly depends on the data normalization method to satisfy positivity and mass balance. We will analyze the dilemma and introduce some normalization methods that are robust in realistic large-scale inversions as well as preserving the convexity. 
%
%In the following sections, we will briefly introduce full waveform inversion (FWI), optimal transport and the merit of using optimal transport based ideas to tackle the current limitations in FWI. Data normalization is the most important step of applying optimal transport to FWI, which will be the primary focus of this paper.

\section{Full Waveform Inversion}
% 2. Full Waveform Inversion (formulation, adjoint state method and problem with cycle skipping)
Full Waveform Inversion (FWI) is a nonlinear inverse technique that utilizes the entire wavefield information to estimate the earth properties. 
The notion of FWI was first brought up three decades ago~\cite{lailly1983seismic, tarantola1982generalized} and has been actively studied and applied with the increase in computing power. It is now a common technique in practice. 

Wave-propagation modeling is the most basic step in seismic imaging. Without loss of generality, we will explain everything in a simple acoustic setting in this paper:
\begin{equation}\label{eq:FWD}
     \left\{
     \begin{array}{rl}
     & m(\mathbf{x})\frac{\partial^2 u(\mathbf{x},t)}{\partial t^2}- \Laplace u(\mathbf{x},t) = s(\mathbf{x},t)\\
    & u(\mathbf{x}, 0 ) = 0                \\
    & \frac{\partial u}{\partial t}(\mathbf{x}, 0 ) = 0    \\
     \end{array} \right.
\end{equation}
We assume the model $m(\mathbf{x}) = \frac{1}{c(\mathbf{x})^2}$ where $c(\mathbf{x})$ is the velocity, $u(\mathbf{x},t)$ is the wavefield, $s(\mathbf{x},t)$ is the source. It is a linear PDE but a nonlinear operator from model domain $m(\mathbf{x})$ to data domain $u(\mathbf{x},t)$. 

As we will see, the mathematical formulation of FWI is PDE constrained optimization. The objective function is the misfit between the synthetic data which is generated by solving certain wave equation numerically with predicted model parameters and the observed data measured from the field which is a result of natural propagation with the real physics. 
%In both time~\cite{tarantola1987inverse} and frequency domain~\cite{PRATT1990,pratt1990inverse}, least-squares norm has been the most widely used misfit function. 
For example, in time domain conventional FWI defines a least-squares waveform misfit as
\begin{equation}
d(f,g) = J_1(m)=\frac{1}{2}\sum_r\int\abs{f(\mathbf{x_r},t;m)-g(\mathbf{x_r},t)}^2dt,
\end{equation}
where $\mathbf{x_r}$ are receiver locations, $g$ is observed data, and $f$ is simulated data which solves~\eqref{eq:FWD} with model parameter $m$. This formulation can also be extended to the case with multiple sources.

In large-scale realistic 3D FWI, there are typically millions of variables describing $m(\mathbf{x})$. 
It is not practical to compute the derivative of the misfit function with respect to each model variable directly.
%is possible with the advances in numerical methods and computational power to solve the 3D wave equation and compute the Fr\'{e}chet derivative with respect to model parameters efficiently.  
With the adjoint-state method, one only needs to solve two wave equations numerically to compute the Fr\'{e}chet derivative, the forward propagation and the adjoint wavefield propagation. Different misfit functions $J(m)$ typically only affect the source term in the adjoint wave equation~\cite{Plessix,tarantola2005inverse}. The gradient is similar to the usual imaging condition~\cite{Claerbout1971}:
\bq~\label{eq:adj_grad3}
\frac{\partial J}{\partial m}  =- \int_0^T \frac{\partial^2 u(\mathbf{x},t)}{\partial t^2} v(\mathbf{x},t)dt,
\eq
where $v$ is the solution to the adjoint wave equation:
\begin{equation} \label{eq:FWI_adj}
     \left\{
     \begin{array}{rl}
     & m\frac{\partial^2 v(\mathbf{x},t)}{\partial t^2}- \Laplace v(\mathbf{x},t)  = R^T\frac{\partial J}{\partial f}\\
    & v(\mathbf{x}, T) = 0                \\
    & v_t(\mathbf{x}, T ) = 0                \\
     \end{array} \right.
\end{equation}
Here $R$ is a restriction operator only at the receiver locations.

\begin{figure}
	\subfloat[Two signals]{\includegraphics[width=1.0\textwidth]{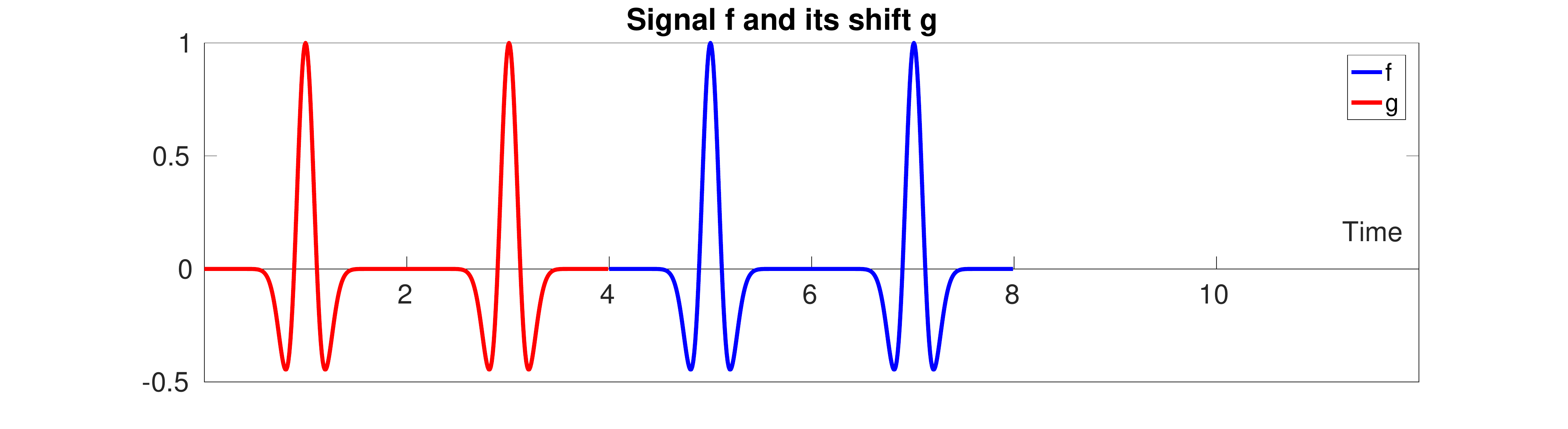}\label{fig:2_ricker_signal}}\\
  	\subfloat[$L^2$ sensitivity curve]
{\includegraphics[width=0.5\textwidth]{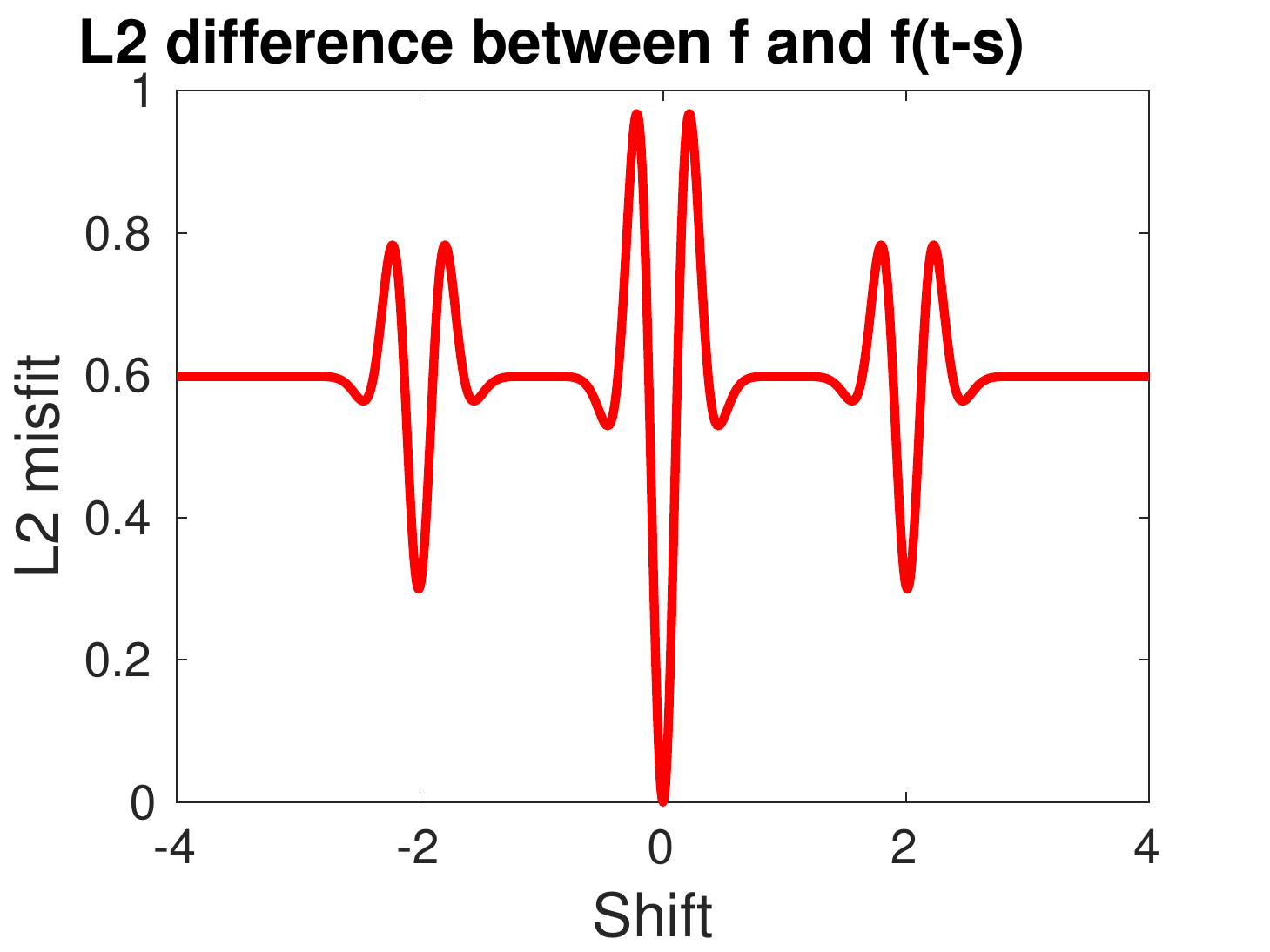}\label{fig:2_ricker_L2}}
  	\subfloat[$W_2$ sensitivity curve]
{\includegraphics[width=0.5\textwidth]{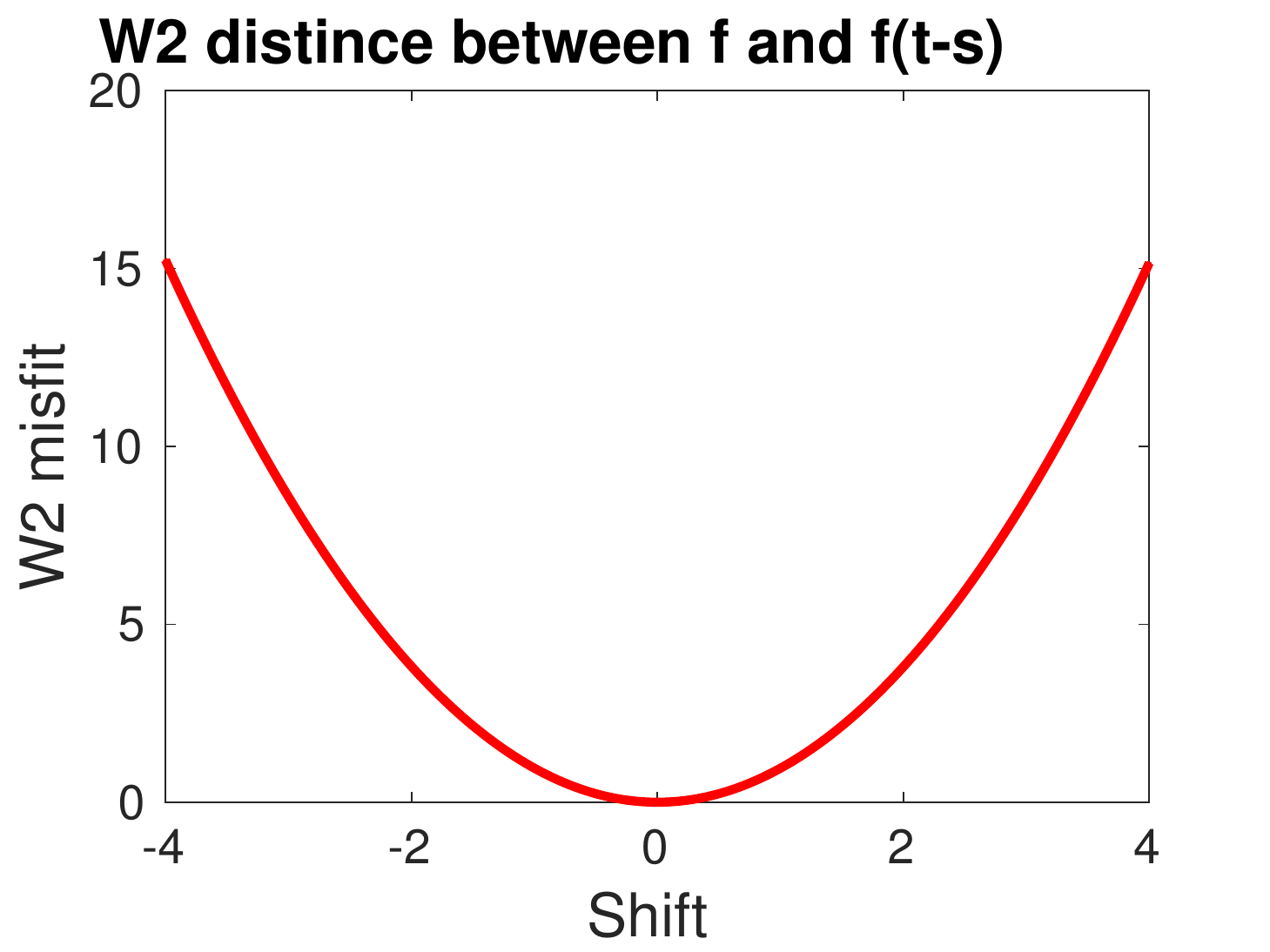}\label{fig:2_ricker_W2}}

\caption{(A)~A signal consisting two Ricker wavelets (blue) and its shift (red) (B)~$L^2$ norm between $f$ and $g$ which is a shift of $f$.  (C) $W_2$ norm between $f$ and $g$ in terms of different shift $s$ }
\end{figure}

It is well known that the accuracy of FWI with $L^2$ norm as misfit functional deteriorates from the lack of low frequencies, data noise, and poor starting model, which may result in local minima trapping. These limitations are on top of the potential ill-posedness of the inverse problem which we here treat as a PDE-constrained optimization.
Figure~\ref{fig:2_ricker_signal} displays two signals, each of which contains two Ricker wavelets and $f$ is simply a shift of $g$. The $L^2$ norm between $f$ and $g$ is plotted in Figure~\ref{fig:2_ricker_L2} as a function of the shift $s$. We observe many local minima and maxima in this simple two-event setting which again demonstrated the difficulty of the so-called cycle-skipping issues~\cite{yang2017analysis}. 

%The lower frequency components have a wider basin of attraction with the least-squares norm being the misfit function. Several hierarchical methods that invert from low frequencies to higher frequencies have been proposed in the literature to mitigate the cycle-skipping of the inverse problem~\cite{bunks1995multiscale, kolb1986pre, pratt1990inverse,sirgue2004efficient,weglein2003inverse}.
A recently introduced class of misfit functions to tackle the cycle-skipping issue is the quadratic Wasserstein metric~\cite{chen2017quadratic,EFWass, engquist2016optimal,yang2017analysis,yangletter, yang2017application}. The $L^2$ misfit function measures the difference in amplitude locally. The optimal transport based methods compare the observed and simulated data globally and thus more effectively include phase information. 

As a useful tool from the theory of optimal transport, the quadratic Wasserstein metric ($W_2$) computes the minimal cost of rearranging one distribution into another with a quadratic cost function. The squared Wasserstein metric has several properties that make it attractive as a choice for misfit function~\cite{engquist2016optimal}.  One highly desirable feature is its convexity with respect to several parameterizations that occur naturally in seismic waveform inversion. As seen in  Figure~\ref{fig:2_ricker_W2}, $W_2$ norm significantly improves the convexity of the misfit sensitivity curve.  
Another important property of optimal transport is the insensitivity to noise. One can find the more theoretical results in~\cite{engquist2016optimal} and the numerical examples in ~\cite{yang2017application}.

\section{Optimal Transport and the Wasserstein metric}
% 3. Wasserstein metric and Optimal transport (definition, some of our theorems, briefly Monge-Ampere, some results, something on unbalanced?)
The topic of optimal transport starts with the problem brought up by Gaspard Monge in 1781~\cite{Monge}.  Let X and Y be two metric spaces with probability measures $\mu$ and $\nu$ respectively. Assume X and Y have equal total measure:
\bq
\int_X d\mu = \int_Y d\nu
\eq
Without loss of generality, we will hereafter assume the total measure to be one, i.e., $\mu$ and $\nu$ are probability measures.

\begin{definition}[Mass-preserving map]
A map $T: X \rightarrow Y$ is mass-preserving if for any measurable
set $B \in  Y$ ,
\bq
\mu (T^{-1}(B)) = \nu(B)
\eq
If this condition is satisfied, $\nu$ is said to be the push-forward of $\mu$ by $T$, and we write $\nu = T_\# \mu $
\end{definition}

Given two nonnegative densities $f = d\mu$ and $g=d\nu$, we are interested in the mass-preserving map $T$ such that $f = g \circ T$. The transport cost function $c(x,y)$ maps pairs $(x,y) \in X\times Y$ to $\mathbb{R}\cup \{+\infty\}$, which denotes the cost of transporting one unit mass from location $x$ to $y$. The most common choices of $c(x,y)$ include $|x-y|$ and $|x-y|^2$. We are interested in finding the optimal map that minimizes the total cost which formally defines a class of metrics: the Wasserstein distance:
\begin{definition}[The Wasserstein distance]
  %The mathematical definition of the Wasserstein distance between the distributions
  We denote by $\mathscr{P}_p(X)$ the set of probability measures with finite moments of order $p$. For all $p \in [1, \infty)$,   
\bq~\label{eq:static}
W_p(\mu,\nu)=\left( \inf _{T_{\mu,\nu}\in \mathcal{M}}\int_{\mathbb{R}^n}\left|x-T_{\mu,\nu}(x)\right|^p d\mu(x)\right) ^{\frac{1}{p}},\quad \mu, \nu \in \mathscr{P}_p(X).
\eq
$\mathcal{M}$ is the set of all maps that rearrange the distribution $\mu$ into $\nu$.
\end{definition}

The optimal transport in higher dimension has no explicit solutions. It is an infinite dimensional optimization problem if we search directly in the function space for $T$. An alternative is to solve the relaxed dual problem by outstanding techniques in linear programming, for example, the alternating direction method of multipliers (ADMM), see the survey~\cite{glowinski2014alternating} by Glowinski. However, the optimal map takes on additional structure in the special case of a quadratic cost function (i.e. $c(x,y) = |x-y|^2$).
The following Brenier's theorem~\cite{Brenier, DePhilippis2013} gives an elegant result about the uniqueness of optimal transport map for the quadratic cost as well as its intrinsic connection with the \MA equation:
\begin{theorem}[Brenier's theorem~\cite{Villani}]
Let $\mu$ and $\nu$ be two compactly supported probability measures on $\R^n$. If $\mu$ is absolutely continuous with respect to the Lebesgue measure, then
\begin{enumerate}
\item There is a unique optimal map $T$ for the cost function $c(x,y) = |x-y|^2$. \item There is a convex function $u: \R^n \rightarrow \R$ such that the optimal map $T$ is given by $T(x) = \nabla u(x)$ for $\mu$-a.e. x.
\end{enumerate}
Furthermore, if $\mu(dx) = f(x)dx$, $\nu(dy) = g(y)dy$, then $T$ is differential $\mu$-a.e. and 
\bq
\det (\nabla T(x)) = \frac{f(x)}{g(T(x))}.
\eq
\end{theorem}

According to Brenier's theorem, in order to compute the misfit between distributions $f$ and $g$, one can first get the optimal map $T(x) = \nabla u(x)$ via the solution of the following \MA equation:
\bq\label{eq:MAA}
\det (D^2 u(x)) = \frac{f(x)}{g(\nabla u(x))}, \quad u \text{ is convex}.
\eq
Typically it is coupled to the non-homogeneous Neumann boundary condition 
\bq\label{eq:BC}
\nabla u(x) \cdot \nu = x\cdot \nu, \,\, x \in \partial X.
\eq
The squared Wasserstein metric is then given by
\bq\label{eq:WassMA}
W_2^2(f,g) = \int_X f(x)\abs{x-\nabla u(x)}^2\,dx.
\eq

We have followed~\cite{benamou2014numerical} for the numerical solution to the \MA equation when computing the quadratic Wasserstein distance for the global comparison in FWI~\cite{yang2017application}. For a survey of recent numerical methods for nonlinear second order PDEs, see~\cite{feng2013recent}.

\section{Data Normalization}
The primary constraints for applying optimal transport to general signals are that the functions should be restricted to nonnegative measures sharing equal total mass (e.g., probability distributions). This is a crucial limitation for many applications that need to compare general signals or allow for only partial displacement of the mass.

\subsection{Background}
There are many proposals in the literature for dealing with the mass balance constraint. Two notions particularly stand out, which are derived rigorously as an extension based on the original optimal transport problem. One is the unbalanced optimal transport, which is formulated as another well-defined metric named the Wasserstein-Fisher-Rao distance~\cite{chizat2016interpolating,kondratyev2016new}. The other approach is the optimal partial transport whose mathematical properties are discussed in detail by~\cite{caffarelli2010free,figalli2010optimal}. As a comparison, there are very few papers discussing the positivity constraint. In~\cite{mainini2012description}, a proposal is made to recombine the data using the decomposition in positive and negative part to compare positive measures with mass conservation.
It is based on the following special dual form of the $W_1$ metric, i.e., $p=1$ in ~\eqref{eq:static}, between  density functions $f = d\mu$ and $g = d\nu$ :
\bq \label{eq:W1}
W_1(f,g)  % = \inf\limits_{T\in\M} \int\limits_X\abs{x-T(x)}f(x)\,dx 
=\max_{\varphi \in \text{Lip}_1} \int_X \varphi(x)(f(x) - g(x))dx, 
\eq
where $\text{Lip}_1$ is the space of all 1-Lipschitz functions.

Based on the dual formulation above, one can easily extend it to signed measures $f$ and $g$ by defining
\begin{eqnarray}
\widetilde{W_1}(f,g)  &=& \widetilde{W_1}(f^+ - f^-, \ g^+ - g^-)\\
&=& \max_{\varphi \in \text{Lip}_1} \int_X \varphi(x)(f^+ - f^- - g^+ + g^-)dx \\
&=& \max_{\varphi \in \text{Lip}_1} \int_X \varphi(x)(f^+  + g^-  -( f^- + g^+ ))dx \\
 &= & \widetilde{W_1}(f^+ + g^- , \ f^- +g^+)\\
 & = & W_1(\rho_1,\rho_2) \label{eq: realW1},
\end{eqnarray}

where $\rho_1 =  f^+ + g^-$, the sum of the positive part of $f$ and the negative part of $g$, and $\rho_2 =  f^- + g^+$, the sum of the negative part of $f$ and the positive part of $g$. The $W_1$ in~\eqref{eq: realW1} is same as the standard 1-Wasserstein distance in ~\eqref{eq:W1}.

\begin{figure}
\includegraphics[width=1.0\textwidth]{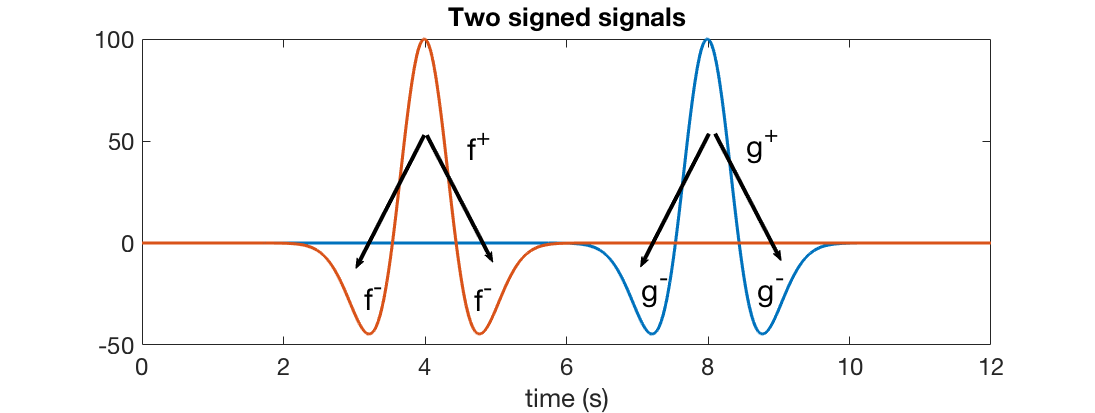}
\caption{The optimal transport may map $f^+$ to $f^-$ and $g^+$ to $g^-$ if formulated as ~\eqref{eq: realW1} (arrows indicate transport) }\label{fig:wrongW1}
\end{figure}

The formulation above defines a cost for transporting signed measures. However, it is not a  canonical optimal transport distance. There is a risk that the true optimal transport represented in~\eqref{eq: realW1} matches $f^+$ to $f^-$ and $g^+$ to $g^-$ under certain circumstances (Figure~\ref{fig:wrongW1}). Especially in FWI, we want to map one signal to the other instead of compensating within one signal itself.

\subsection{Early ideas}
In this section, we will introduce some normalization ideas which we proposed in the past to transform seismic data into probability signals such that the standard optimal transport theory will apply. We will analyze their properties and in particular why they often have problems with realistic large-scale FWI.

\begin{figure}
\includegraphics[width=1.0\textwidth]{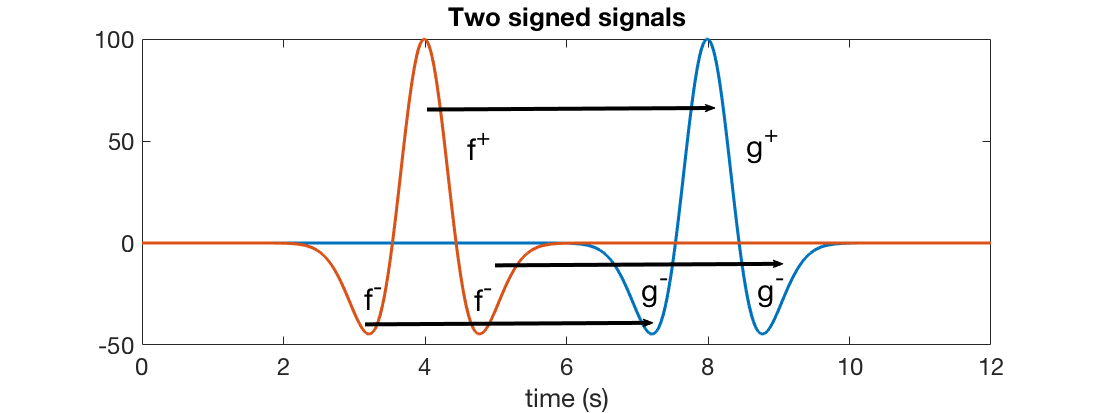}
\caption{The optimal transport plan maps ${f}^+$ to ${g}^+$ and ${f}^-$ to ${g}^-$ if formulated as ~\eqref{eq:separate} (arrows indicate transport) }\label{fig:separate}
\end{figure}

In~\cite{EFWass,engquist2016optimal}, the signals were separated into positive and negative parts $f^+ = \max\{f,0\}$, $f^- = \max\{-f,0\}$ and scaled by the total mass $\langle f \rangle = \int_X f(x)\,dx$ (Figure~\ref{fig:separate}).  Inversion was accomplished using the modified misfit function
\bq \label{eq:separate}
J_2(m) = W_2^2\left(\frac{f^+}{\langle f^+ \rangle},\frac{g^+}{\langle g^+ \rangle} \right) + W_2^2\left(\frac{f^-}{\langle f^- \rangle}, \frac{g^-}{\langle g^- \rangle}\right). 
\eq

Recall the adjoint-state equation introduced earlier~\eqref{eq:FWI_adj}. In order to compute the gradient of the misfit function with respect to the model parameters for FWI, we simply need the Fr\'{e}chet derivative of the misfit function with respect to the synthetic data $f$. Therefore, the critical element in the backpropagation is $\frac{\partial J}{\partial f}$. Once we separate the signals as in ~\eqref{eq:separate}, discontinuities are introduced in derivatives of $f$ which causes problems in the optimization process and for the wave equation solvers.
The same principle applies to the absolute-value scaling $W_2^2(\abs{f}, \abs{g})$ since  absolute-value function is not differentiable at zero.  

\begin{figure}
\includegraphics[width=1.0\textwidth]{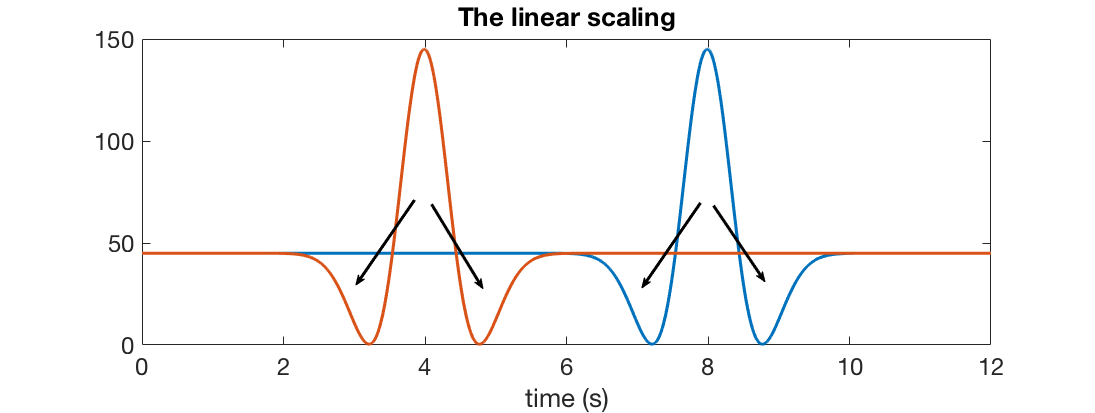}
\caption{The linear scaling: $f\rightarrow f+c$ and $g\rightarrow g+c$; there is chance of having local transport (arrows indicate transport)}\label{fig:linear}
\end{figure}

The linear scaling we used in our earlier papers, i.e., Equation~\eqref{eq:linear}, on the other hand, works very well even if the related misfit lacks strict convexity with respect to shifts (see Figure~\ref{fig:linear}). Here are several beneficial properties about the linear scaling. First, it has a wider basin of attraction than $L^2$ norm when it comes to simple shifts~\cite{yangletter}. The two-variable example described in~\cite{yang2017analysis} is based on the linear scaling. It gives the convexity with respect to a subset of model variables in velocity compared to the result of $L^2$. Second, it provides a smooth bijection between the original data and the normalized data, which is favorable when combining with the adjoint-state method. Third, realistic seismic data always has the mean-zero property after a standard data processing. This indicates that $<f+c>$ is equal to $<g+c>$.
This means that if two short seismic signals or, so-called events, are well matched between f and g they will stay so even after the normalization process and not be influenced by other events further away. The property is essential in the early iteration steps when the simulated signals do not include all details that are in the measured signal. On the other hand, if the individual events are void of zero frequencies the transport defining $W_2$ may be local as is seen in Figure~\ref{fig:linear}, which can cause trapping in local minima.

% 4. Normalizations (this could be the main part, some results, theory and results for "mixed normalization")

\begin{figure}
\includegraphics[width=1.0\textwidth]{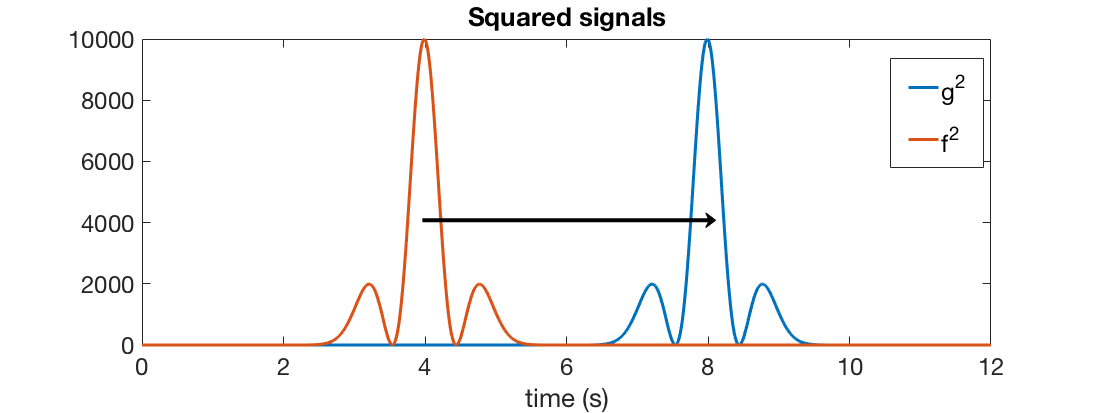}
\caption{Square of the data: $f\rightarrow f^2$ and $g\rightarrow g^2$ (arrows indicate transport)}
\label{fig:f2}
\end{figure}

One scaling method which theoretically should work well with the adjoint-state method is to square the signals first and normalize it to be mass balanced:
\bq\label{eq:f2}
J_3(m) = W_2^2(\frac{f^2}{<f^2>} , \frac{g^2}{<g^2>} ).
\eq
As seen in Figure~\ref{fig:f2}, the two curves are the squares of the two functions in Figure~\ref{fig:separate}. This particular normalization keeps the convexity of the quadratic Wasserstein metric concerning simple shifts like the setting in Figure~\ref{fig:2_ricker_W2}. In~\cite{chen2017quadratic},  squaring the data was used as the normalization to recover a four-variable linear source inversion. However, it has been puzzling since this normalization rarely works well in large-scale inversions with thousands of variables, such as the Camembert example and the standard Marmousi benchmark which we will show later.

\begin{figure}
\centering
  \subfloat[]{\includegraphics[width=0.5\textwidth]{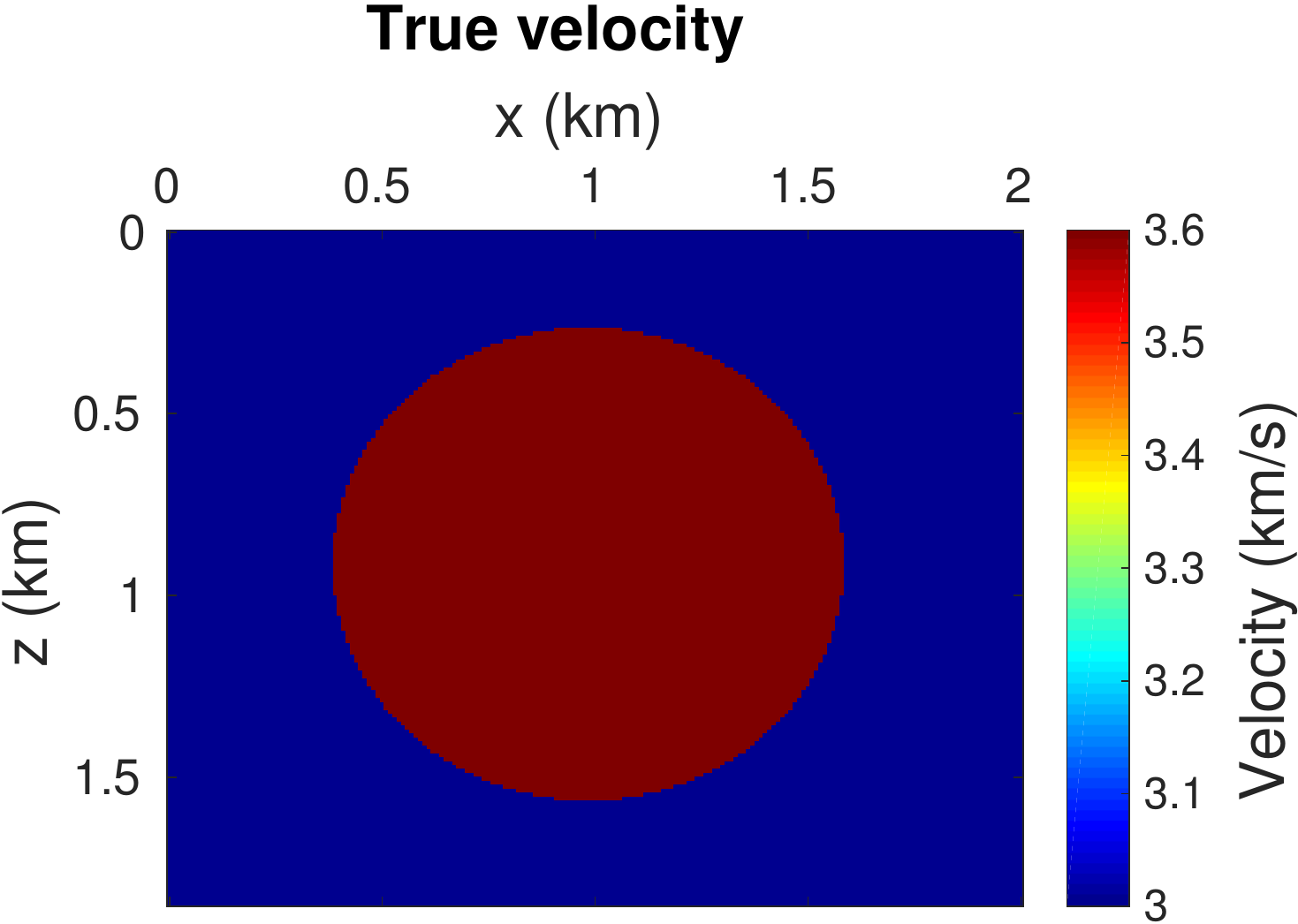}\label{fig:cheese_true}}
  \subfloat[]{\includegraphics[width=0.5\textwidth]{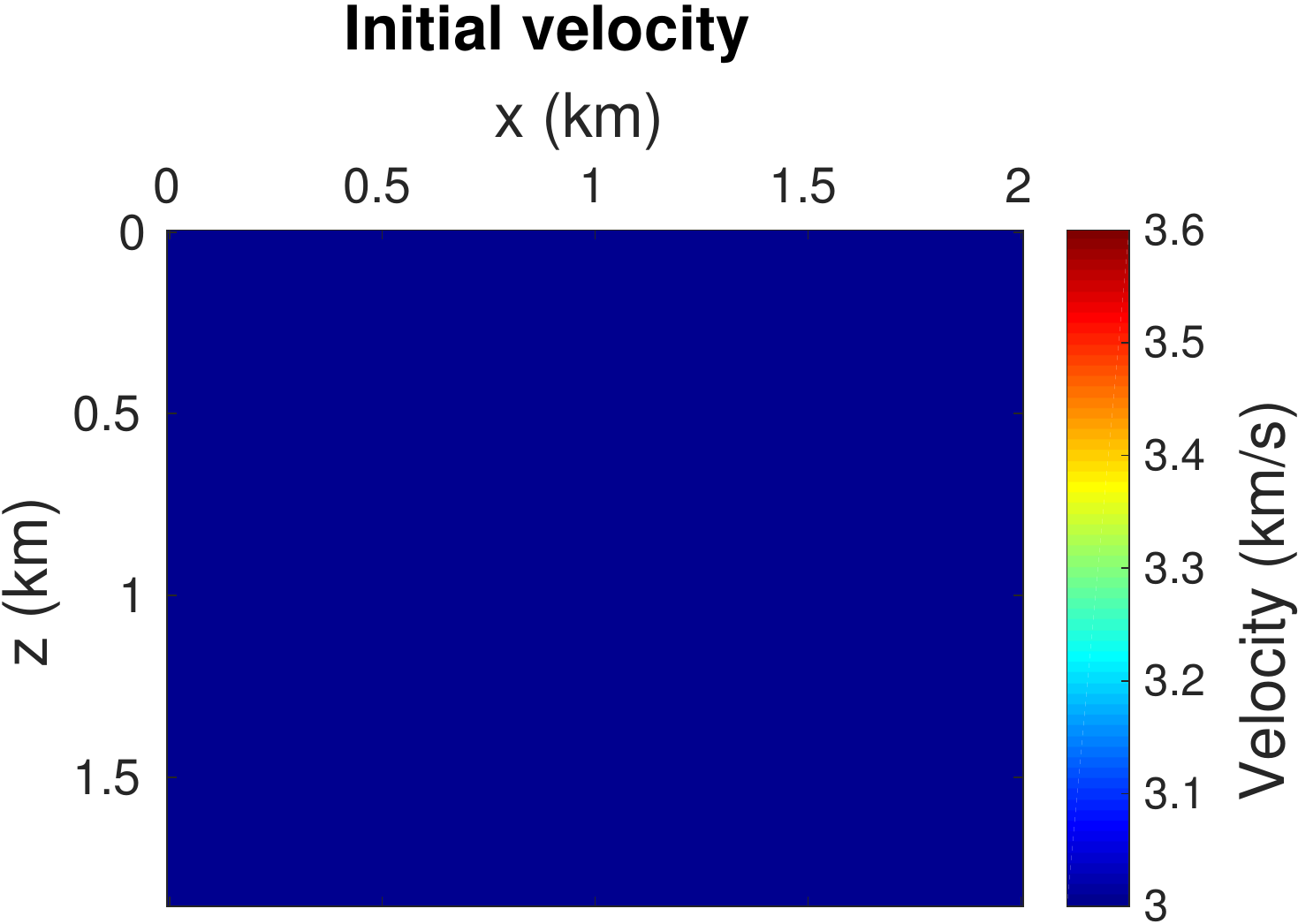}\label{fig:cheese_v0}}
  \caption{(a)~True velocity and (b)~inital velocity for the Camembert model}
  \label{fig:cheese_true,cheese_v0}
\end{figure}

\begin{figure}
      \subfloat[Gradient  by $\tilde{f} = af+b$]
{\includegraphics[width=0.5\textwidth]{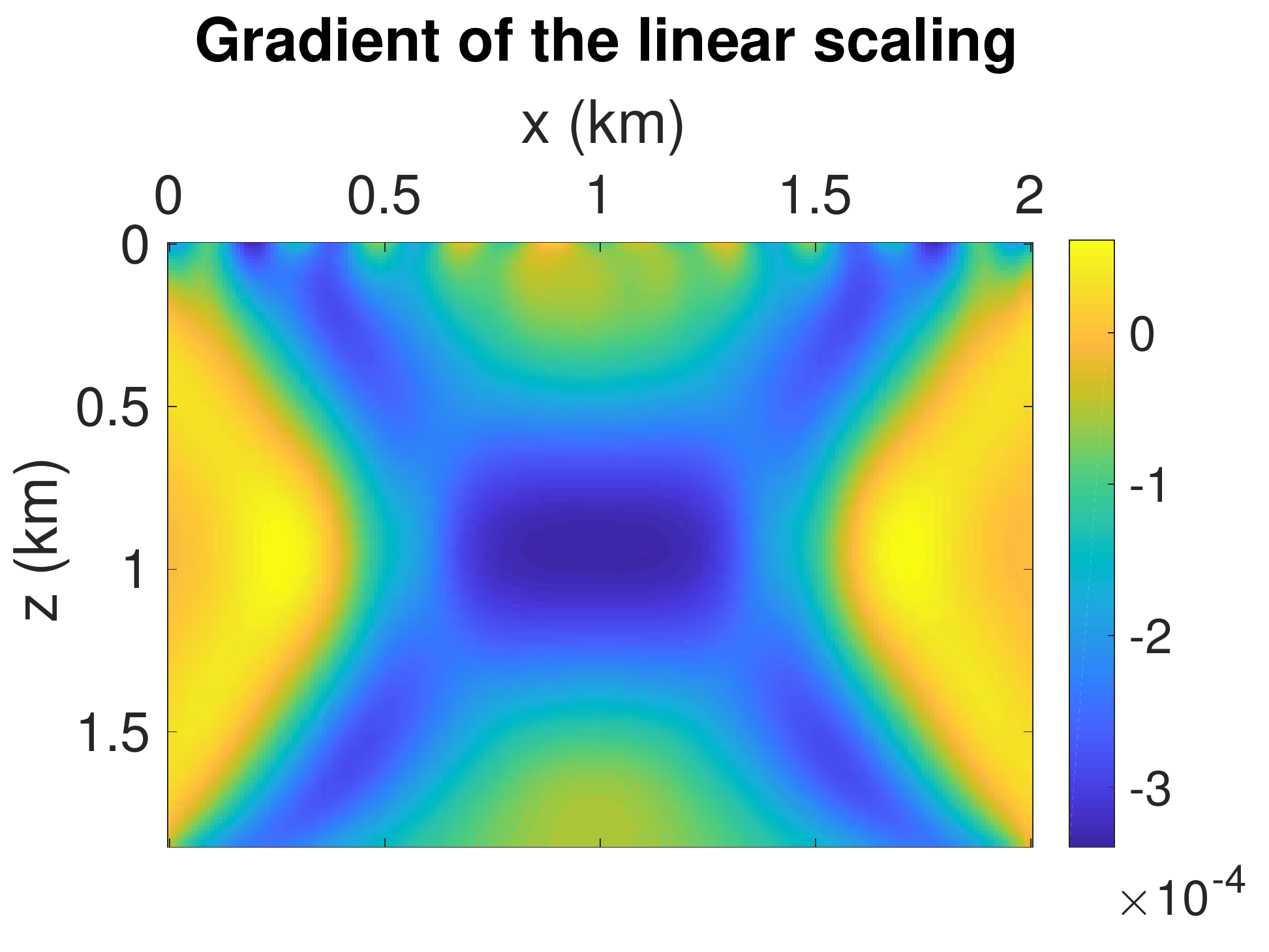}\label{fig:lin_grad}}
      \subfloat[Gradient  by $\tilde{f} = f^2$]
{\includegraphics[width=0.5\textwidth]{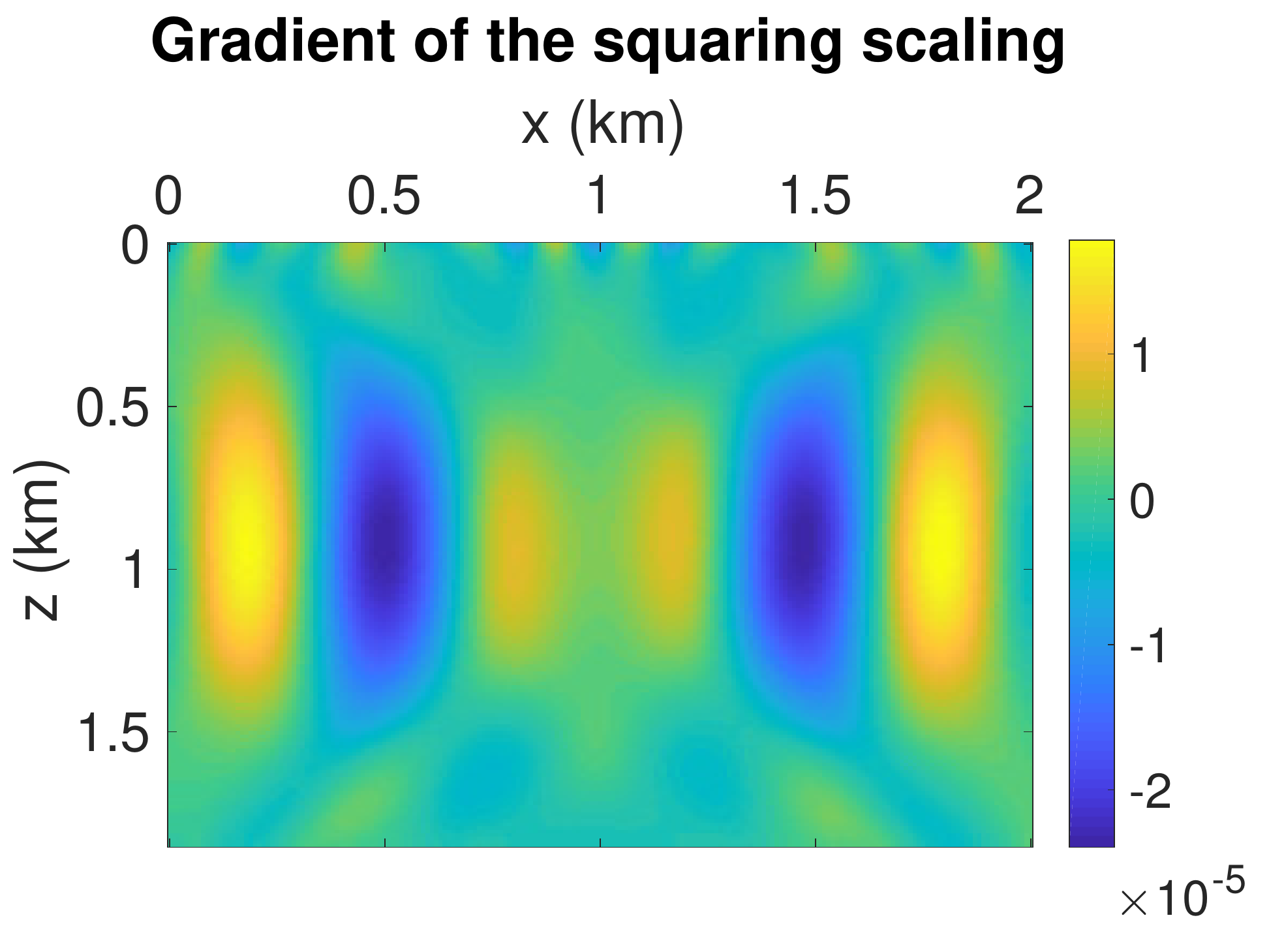}\label{fig:f2_grad}}

\caption{The gradient in the first iteration of the inversion by using~ (a) the linear scaling as the data normalization, and (b) the squaring scaling as the data normalization.}~\label{fig:cheese_grad}
\end{figure}

\begin{figure}
      \subfloat[Inversion  by $\tilde{f} = af+b$]
{\includegraphics[width=0.5\textwidth]{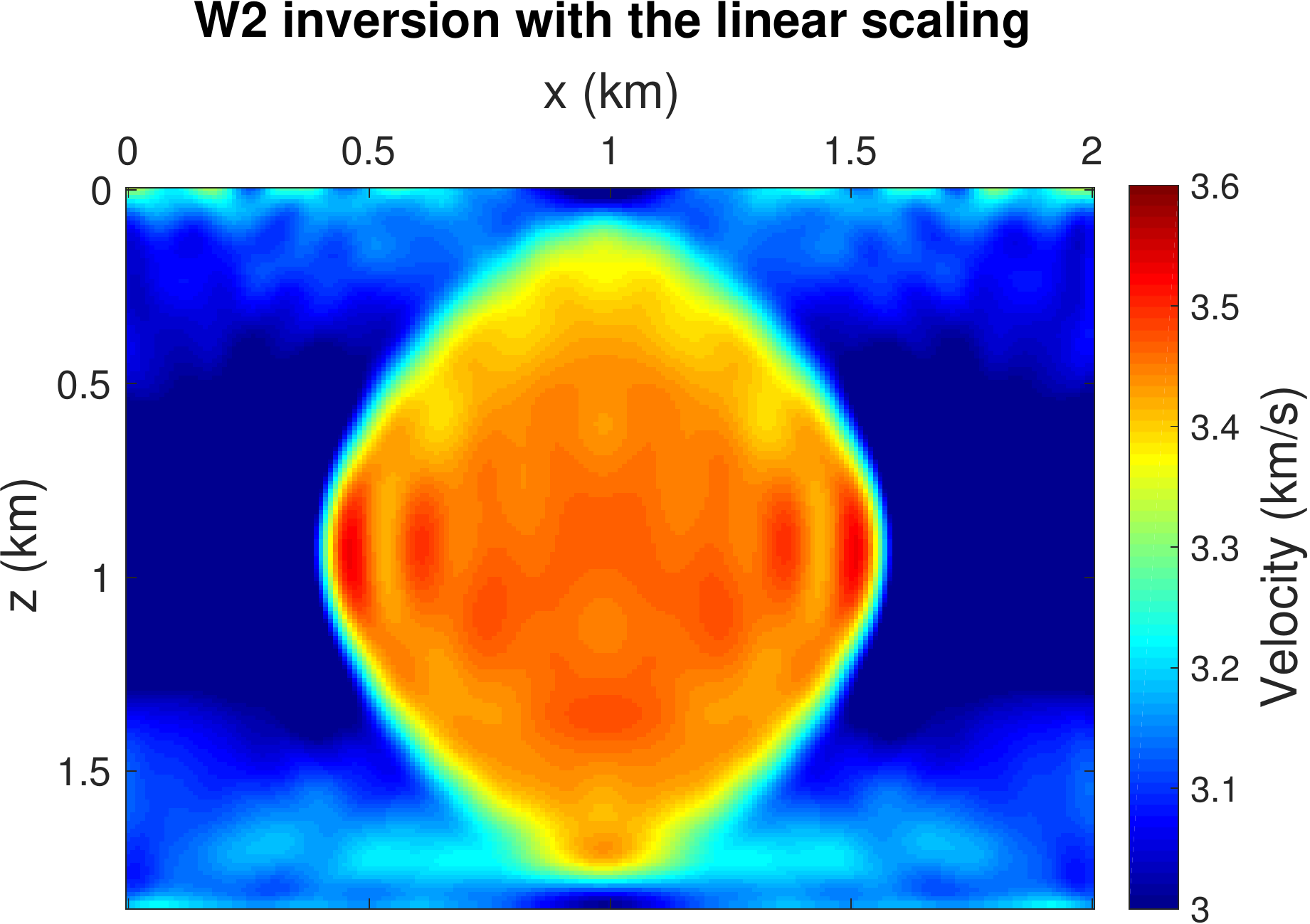}\label{fig:lin_inv}}
      \subfloat[Inversion  by $\tilde{f} = f^2$]
{\includegraphics[width=0.5\textwidth]{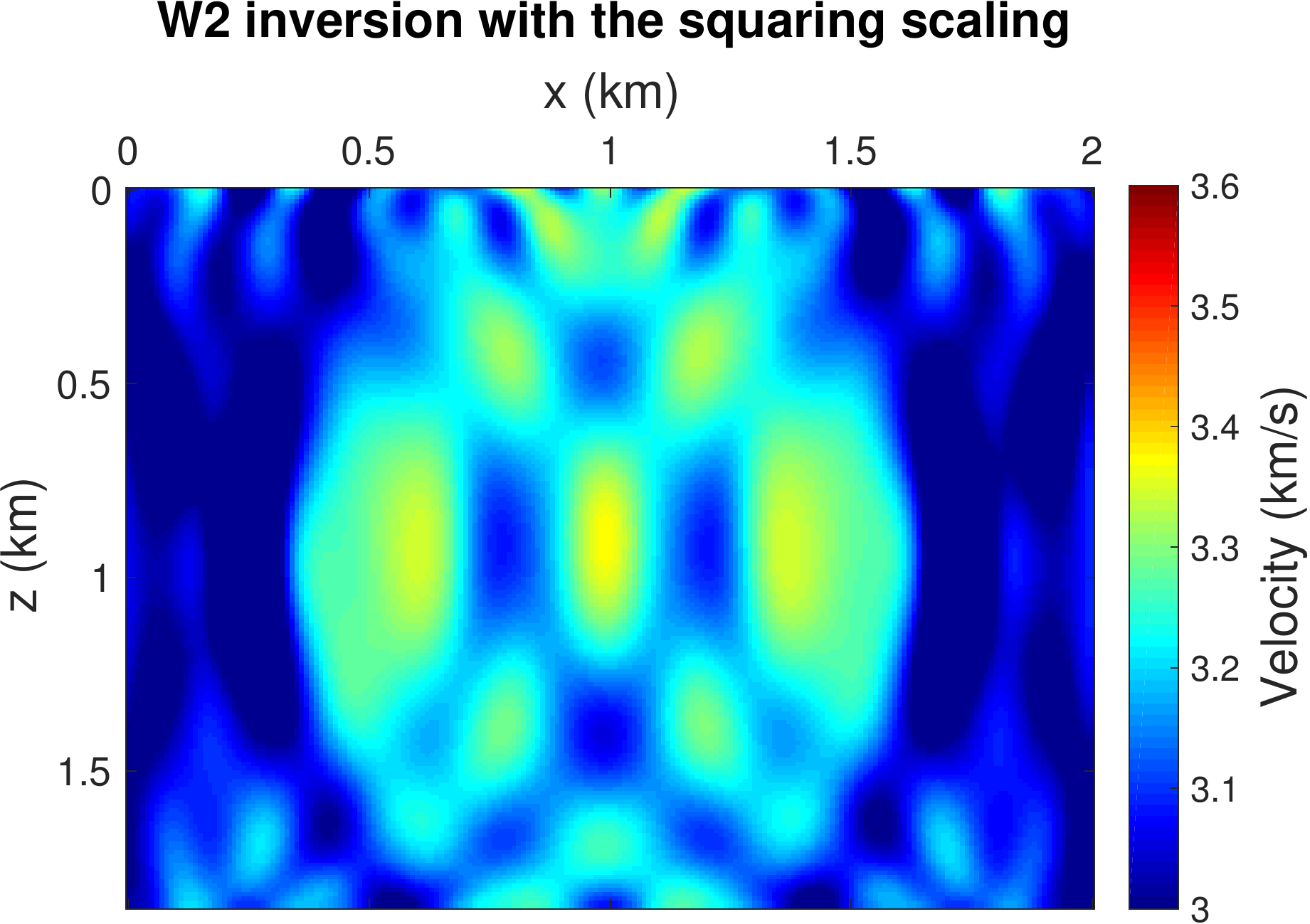}\label{fig:f2_inv}}
\caption{(a)~The inversion result of using ~\eqref{eq:linear} as the data normalization (b)~The inversion result of using ~\eqref{eq:f2} as the data normalization. }\label{fig:cheese_inv}
\end{figure}

It has been a dilemma until recently we can point to three potential factors that may lead to the difficulties. First of all, taking the squares boosts the higher frequency of the signal. It is well known that FWI becomes more difficult as the frequency increases. The robust convergence range is typically within half wavelength~\cite{virieux2009overview}. Just consider a simple oscillartory pulse $\sin(t)^2 = (1-\cos(2t))/2$. 
Second, the refracted, or so-called, diving wave and the reflection wave may reach the receiver at the same time with a similar amplitude but  entirely different polarity. The positive and negative parts of the signal are interchanged. Squaring the signals may lose the important phase information here. Third, when we are dealing with one event in $f$ and multiple events in $g$, $<g^2>$ can be significantly larger than $<f^2>$. 
This is often the case in the initial state of inversion when only one or a few reflections interfaces are known. The measured data $g$ naturally contains the effect of all reflections.
The mass normalization step can distort the correct parts of the signals that both $f$ and $g$ share, and consequently leads to a wrong update in the correct model variables.

%The last and the most significant issue is that squaring the data may create more non-uniqueness and consequently more local minima. For example, for any simple domain decomposition $\Omega = \Omega_1 \cup \Omega_2$, there exists four global minima for synthetic data $f$ and consequently four different models $m_i$ ($ i= 1,2,3,4$) that can generate the zero misfit:
%\bq
%     f^\star = \left\{
%     \begin{array}{rl}
%     &  \pm g\quad \text{on\ } \Omega_1   \\
%     & \pm g \text{\quad on\ } \Omega_2 \end{array} \right. 
%\eq
%The objective function $J_3(m)$ in ~\eqref{eq:f2} is the optimal cost of mapping $\frac{f^2}{<f^2>} $ to $\frac{g^2}{<g^2>}$. The gradient $\frac{\partial J}{\partial m}$ represents the fastest way to decrease the objective function. However, it is very possible that the steepest descent direction can converge to any of the  global minima $m_i$'s although only one model $m^\star$ that makes $f^\star = g$, the desirable solution. It is easy to observe more failures in inversion since there are also local minima around each $m_i$.

We want to demonstrate the issues above with a Camembert model. The true velocity is shown in Figure~\ref{fig:cheese_true} and the initial velocity we use in the inversion is Figure~\ref{fig:cheese_v0}. Figure~\ref{fig:cheese_grad} illustrates a comparison in gradients of the first iteration between the linear scaling and the squaring scaling as different normalizations in optimal transport FWI. There are wrong features in Figure~\ref{fig:f2_grad} even in the first iteration. The final inversion results are shown in Figure~\ref{fig:cheese_inv}. The linear scaling converges to a reasonably well model (Figure~\ref{fig:lin_inv}) while squaring the data in normalization leads the inversion to a local minimum (Figure~\ref{fig:f2_inv}). In both of these two experiments, we use the trace-by-trace techqniue (1D optimal transport) to compute the $W_2$ distance between the synthetic data and the observed data.

\subsection{A sign-sensitive normalization}
Based on the analysis of the squaring scaling in the previous section, a bijection between the original data and the normalized data is essential not to deteriorate the ill-posedness of the inverse problem.
An exponential based normalization was proposed in~\cite{qiu2017full} to transform seismic signals to probability distributions:
\bq\label{eq:exp}
\tilde{f}(t) = \frac{\exp(cf(t))}{<\exp(cf)>},\ \tilde{g}(t) = \frac{\exp(cg(t))}{<\exp(cg)>},\ c >0.
\eq
%An appropriate choice of $c$ will keep the convexity of $W_2$ with respect to shifts as the squaring scaling does, but it is a one-to-one mapping that we can recover the original signal from the normalized signal. 

\begin{figure}
\includegraphics[width=1.0\textwidth]{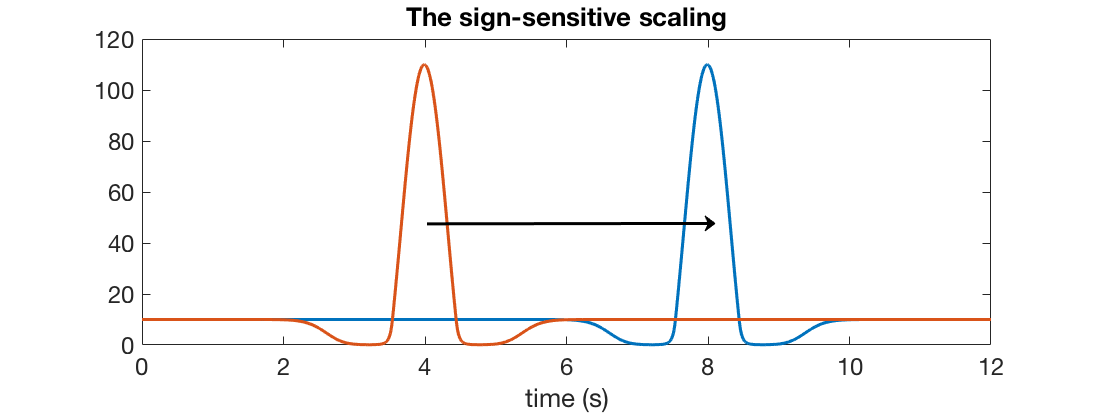}
\caption{One normalization combining both the linear and the exponential methods (Equation~\eqref{eq:mixed})}\label{fig:mixed}
\end{figure}

\begin{figure}
\includegraphics[width=1.0\textwidth]{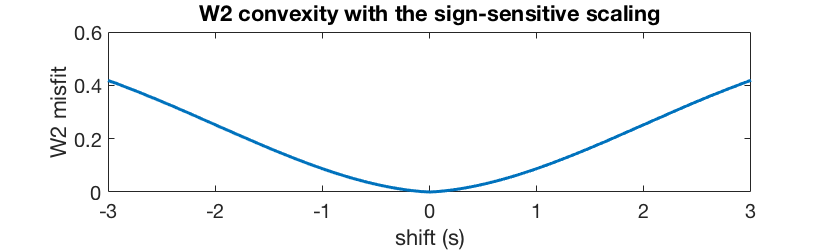}
\caption{The $W_2$ misfit regarding signal shift $s$ by using the sign-sensitive scaling, i.e., $W_2^2(f(t-s),f(t))$.}\label{fig:cvx_mixed}
\end{figure}

Here we propose a new normalization (Equation~\eqref{eq:mixed}) that satisfies most of the essential properties.
%It is a combination between the linear scaling (Equation~\eqref{eq:linear}) and the exponential scaling  (Equation~\eqref{eq:exp}). For the positive part of the signal, we use a linear scaling, and an exponential function is applied to the negative part of the signal~(Figure~\ref{fig:mixed}). 
This normalization ~(Figure~\ref{fig:mixed}) can be seen as a compromise between the linear scaling \eqref{eq:linear} and the positive part scaling ~\eqref{eq:separate}. For the limit of small $c >0$, it is a linear scaling, which directly follows from Taylor expansion of the exponential part. In the limit of large $c$ values, $f$ obviously converges to $f^+$.
It is a $C^1$ function which is compatible with the adjoint-state method. It keeps the convexity of the quadratic Wasserstein distance with respect to signal shifts as shown in Figure~\ref{fig:cvx_mixed}. One needs to select the coefficient $c$ based on the data range. The sign-sensitive scaling is similar to the exponential scaling~\eqref{eq:exp} by suppressing the negative part of the signal, but it does not have the risk of exaggerating large $f$ and $g$ values in the exponential normalization.
%The mixed scaling~\eqref{eq:mixed} used in the quadratic Wasserstein distance can be close to $W_2^2(\frac{f^+}{<f^+>}, \frac{g^+}{<g^+>})$ if $c$ is large, and similar to the linear scaling~\eqref{eq:linear} when $c$ is very small.

If we denote the normalization function in ~\eqref{eq:mixed} as an operator $P$. One idea to use all the information of the signal (instead of just $f^+$) is to consider the following objective function:
\bq\label{eq:mixmix}
J_4(m) = W_2^2 (P(f), P(g)) + W_2^2(P(-f), P(-g)).
\eq

\begin{figure}
\centering
  \subfloat[]{\includegraphics[width=0.45\textwidth]{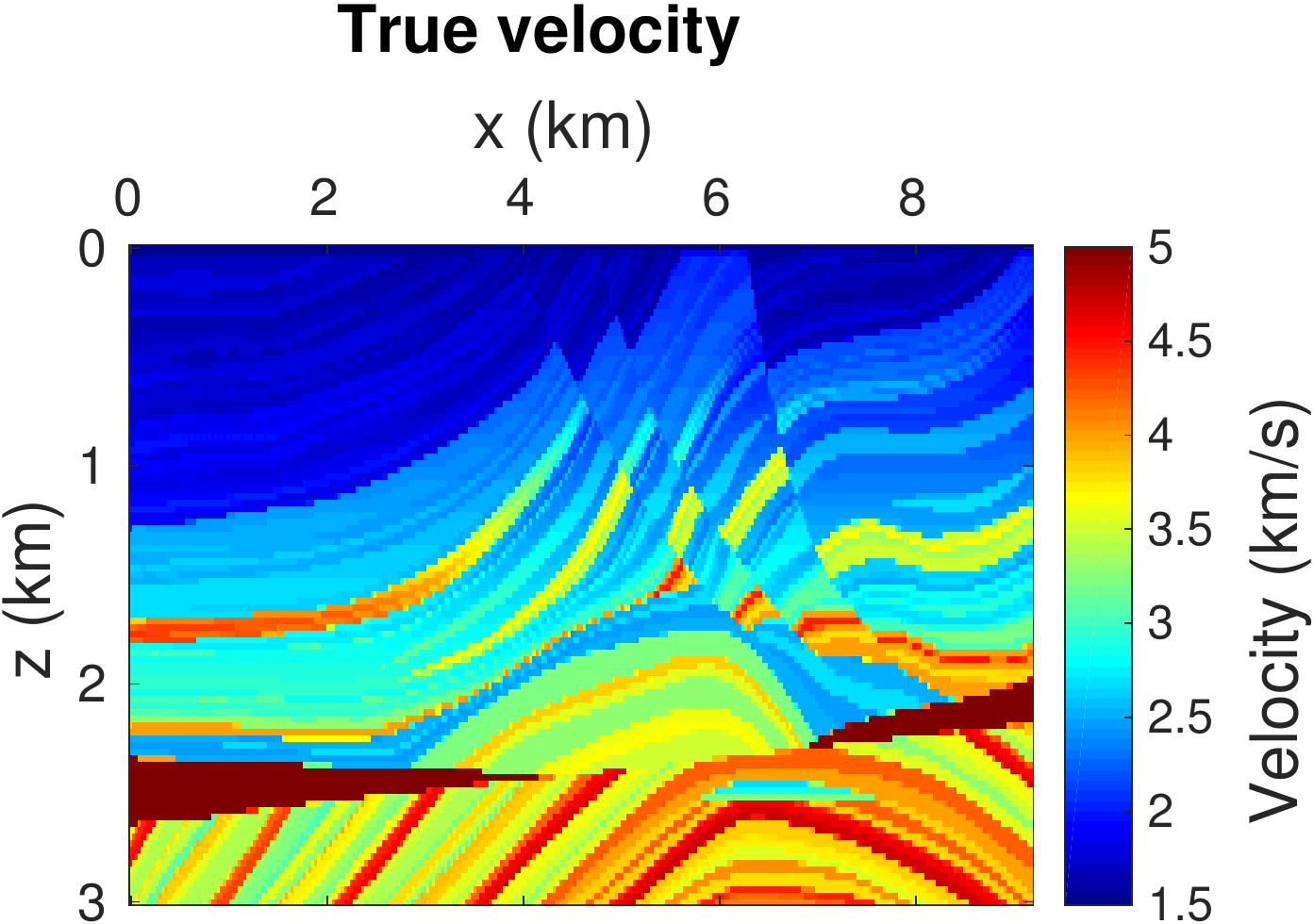}\label{fig:marm2_true}}
  \subfloat[]{\includegraphics[width=0.45\textwidth]{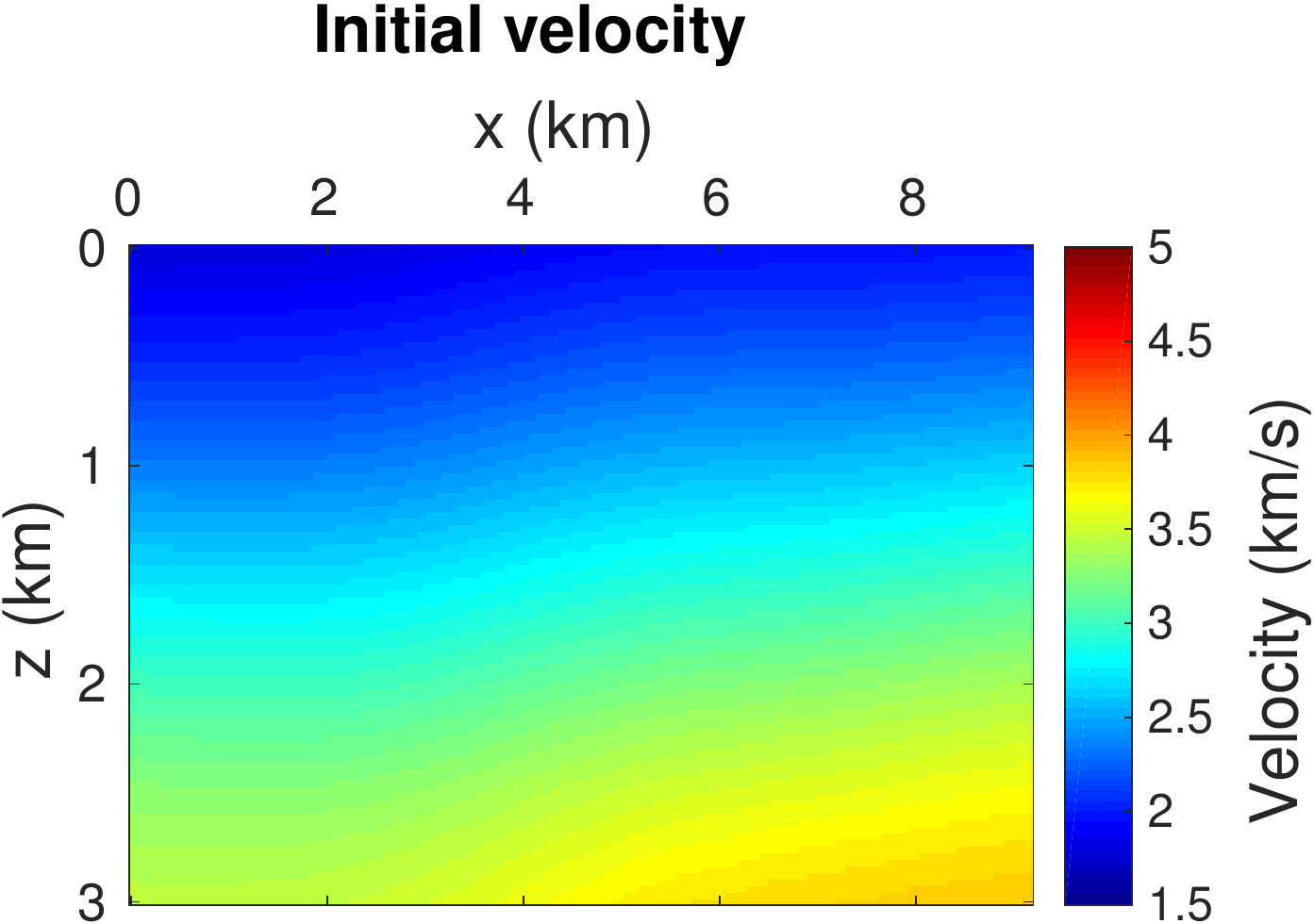}\label{fig:marm2_v0}}
  \caption{(a)~True velocity and (b)~inital velocity for full Marmousi model}
  \label{fig:marm2_true,marm2_v0}
\end{figure}

\begin{figure}
\centering
  \subfloat[]{\includegraphics[width=0.45\textwidth]{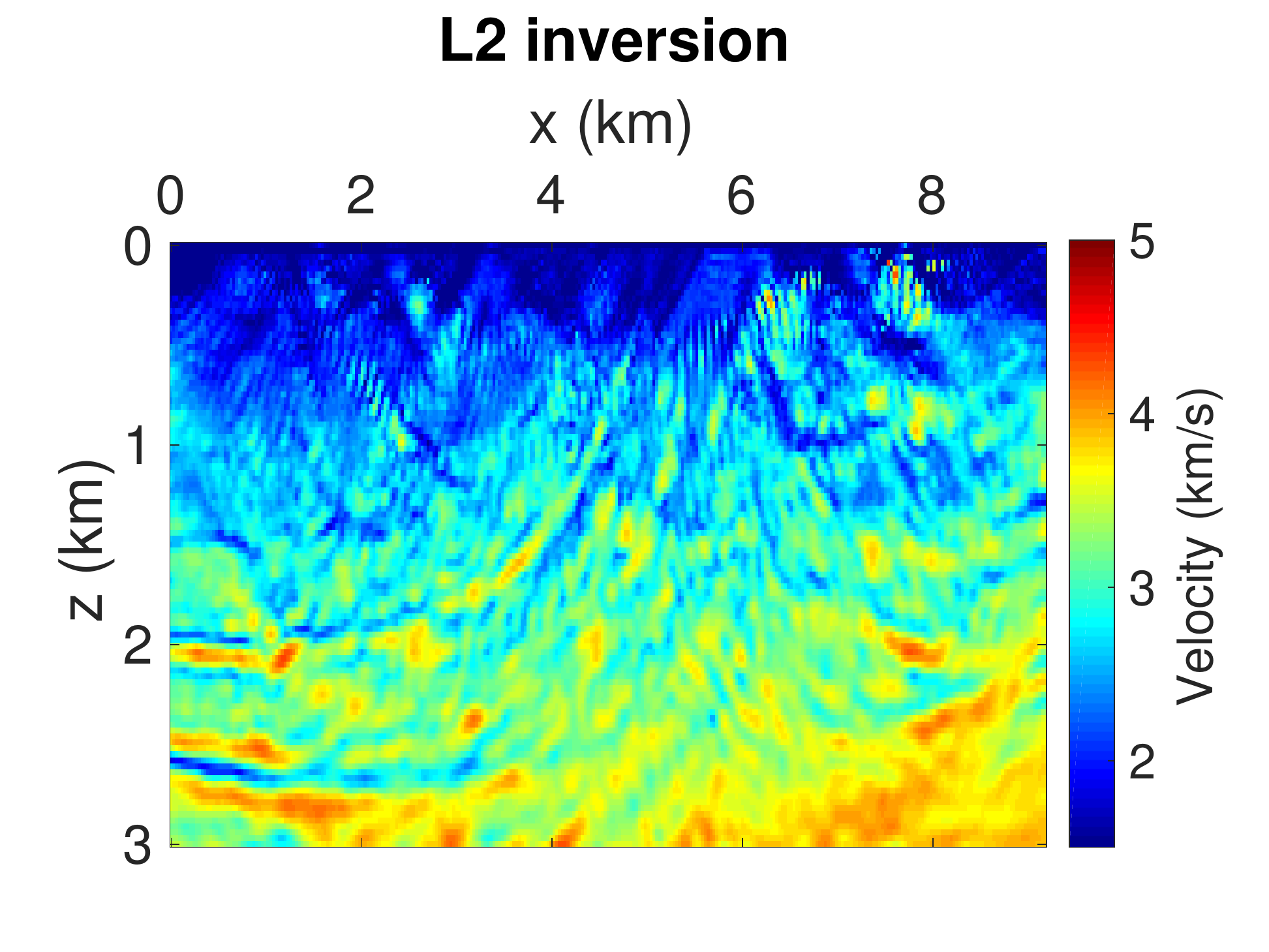}\label{fig:marm2_L2}}
  \subfloat[]{\includegraphics[width=0.45\textwidth]{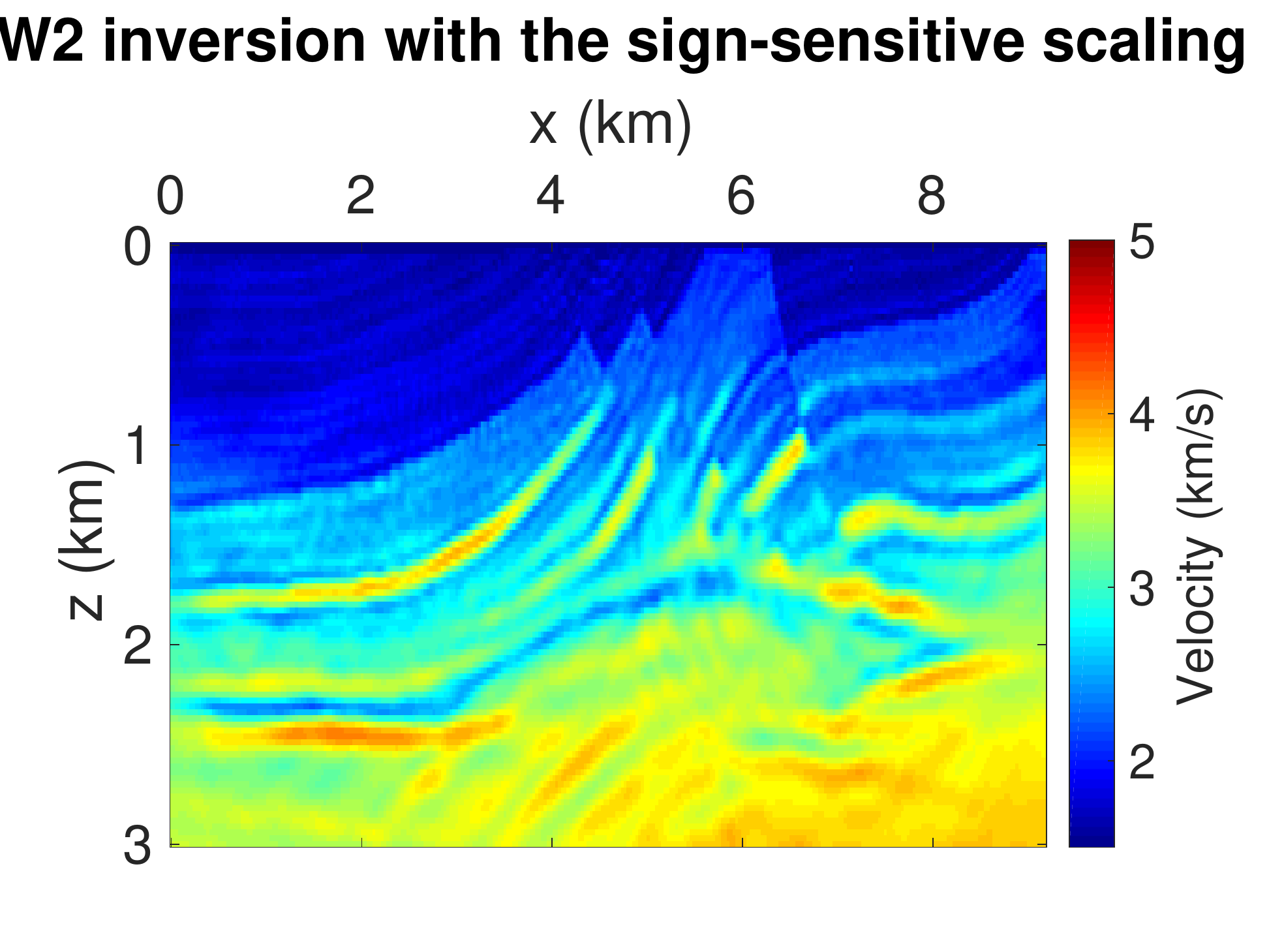}\label{fig:marm2_mixed}}
  \caption{Inversion results of (a)~$L^2$ and (b)~trace-by-trace $W_2$ with the sign-sensitive scaling~\eqref{eq:mixmix}}
  \label{fig:marm2_inv}
\end{figure}

The final experiment is to invert full Marmousi model by conventional $L^2$ and trace-by-trace $W_2$ misfit with $J_4(m)$ as the actual objective function. Figure~\ref{fig:marm2_true} is the P-wave velocity of the full Marmousi model, which is 3km in depth and 9km in width. The inversion starts from an initial model that is the true velocity smoothed by a Gaussian filter with a deviation of 40 (Figure~\ref{fig:marm2_v0}). We place 11 evenly spaced sources on top at 150m depth in the water layer and 307 receivers on top at the same depth with a 30m fixed acquisition. The discretization of the forward wave equation is 30m in the $x$ and $z$ directions and 30ms in time. The source is a Ricker wavelet with a peak frequency of 15Hz, and a high-pass filter is applied to remove the frequency components from 0 to 2Hz. 
Inversions are terminated after 300 l-BFGS iterations which take about 3 hours on a normal workstation. Figure~\ref{fig:marm2_L2} shows the inversion result using the traditional $L^2$ least-squares method and Figure~\ref{fig:marm2_mixed} shows the final result using trace-by-trace $W_2$ misfit function. Again, the result of $L^2$ metric has spurious high-frequency artifacts while $W_2$ using the sign-sensitive scaling~\eqref{eq:mixmix} correctly inverts most details in the true model. 

%\begin{figure}
%\includegraphics[width=1.0\textwidth]{exp.png}\label{fig:exp}
%\caption{The linear scaling: $f\rightarrow \exp(cf)$ and $g\rightarrow \exp(cg)$; there is chance of having local transport.}
%\end{figure}

\section{Conclusion}
Full waveform inversion for seismic imaging and the application of optimal transport for computing the misfit between simulated and measured data are summarized. Seismic signals need to be transformed by some normalization to satisfy the requirements from optimal transport. Advantages and disadvantages of different normalization techniques are discussed. The dilemma that methods, which have provable desirable properties for simple model problems do not work well in practical large-scale settings and other methods that theoretically fail for simple examples perform very well in practice is illuminated. Quadratic scaling belongs to the first class, and linear scaling belongs to the second class of normalizations, which do very well in realistic tests. A new sign-sensitive normalization aiming at bridging these two classes is introduced, and numerical examples are presented.

\bibliographystyle{plain}
\bibliography{Roland}

\end{document}